# 一组 3 个多联立方体平移密铺 3 维空间的不可判定性

杨超[1*], 张诸俊[2]

1. 广东外语外贸大学数学与统计学院, 广州 510006
2. 奉贤区大数据中心, 上海 201499
* 通信作者. E-mail: yangchao@gdufs.edu.cn

**摘要** 平移密铺问题是数学各领域中最基本、最有代表性的不可判定问题之一. Greenfeld 和陶哲轩近几年关于平移密铺的不可判定性有两项重要的工作. 一是证明了在维数充分大的空间中，存在只能非周期地平移密铺空间的单一瓷砖. 二是在空间的维数是问题输入的一部分时，单一瓷砖能否平移密铺空间的周期子集是不可判定的. 这两项工作支撑了如下猜想：存在一个固定的维数 $n$，使得 $n$ 维空间的单一瓷砖平移密铺问题是不可判定的. 解决此猜想的一个途径是对更小的正整数 $k$，证明某固定维数空间中 $k$ 个瓷砖的平移密铺问题的不可判定性, 向 $k=1$ 逼近. 本文证明 3 维空间中一组 3 个多联立方体的平移密铺问题是不可判定的.

**关键词** 不可判定性, 平移密铺, 多联立方体, 王浩瓷砖，Colomb 尺

**MSC (2020) 主题分类** 52C22, 68Q17

## 1 引言

数学各研究领域中都存在若干最基本的问题从可计算理论的角度看是不可判定的，其中组合和离散几何领域中的一个典型代表就是平移密铺问题[10]. 对平移密铺问题的不可判定性的认识源于 1960 年代王浩提出的一个密铺问题[12]. 一个每条边都被赋予某种颜色的正方形称为一个王浩瓷砖，简称王砖. 图1是一个王砖组的例子. 用王砖密铺平面时，要求边对边，且相邻的两个王砖的公共边的颜色要相同. 王浩提出问题：是否存在一个算法，对任意一组王砖，在有限时间内判断能否用这组王砖的平移拷贝（即每一个王砖都可以任意多次平移重复使用），密铺整个平面? Berger 于 1966 年证明了王浩瓷砖平移密铺问题的不可判定性[2].

**定理 1** (见 [2]) 王浩瓷砖的平移密铺问题是不可判定的.

对单一瓷砖的平移密铺问题，研究者长期以来认为问题是可判定的. Lagarias 和汪扬曾明确提出蕴含了可判定性的周期密铺猜想[8]：在任何维数的空间中，单一瓷砖若能密铺整个空间，则总有周期的密铺方案. 出乎意料的是，周期密铺猜想被 Greenfeld 和陶哲轩于 2024 年否证[5]. Greenfeld 和陶哲轩甚至进一步证明了：在一定条件下，单一瓷砖平移密铺空间的不可判定性[6].

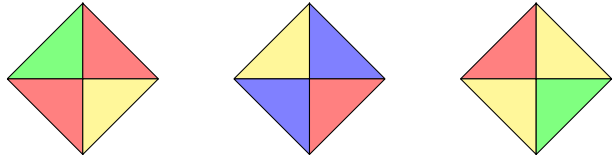

**图 1** 一组 3 个王砖
**Figure 1** A set of 3 Wang tiles

然而，Greenfeld 和陶哲轩的单一瓷砖平移密铺不可判定性[6]有两个额外条件. 一是空间的维数是问题输入的一部分，即空间维数不固定. 二是把问题从密铺整个空间扩充到空间的周期子集. 故以下猜想仍然是一个未解问题.

**猜想 1** 存在一个固定的维数 $n$，使得 $n$ 维空间的单一瓷砖平移密铺问题是不可判定的.

为了回答这一未解问题，一个自然的路径就是寻求更小的正整数 $k$，使得如下固定参数的平移密铺问题是不可判定性的.

**定义 1** ($\mathbb{Z}^n$ 中 $k$ 个瓷砖的平移密铺问题)  设 $n, k$ 为两个固定的正整数，$\mathbb{Z}^n$ 的一个有限子集称为一个瓷砖. 是否存在算法可在有限时间内判定任给的一组 $k$ 个瓷砖能否平移密铺整个 $n$ 维空间 $\mathbb{Z}^n$.

定义1中 $\mathbb{Z}^n$ 的平移密铺问题也可以等价地表述为 $\mathbb{R}^n$ 中的平移密铺问题，其中一个瓷砖为有限个形如 $\Pi_{1 \leqslant i \leqslant n}[z_i, z_i+1], (i \in \mathbb{Z})$ 的单位超立方体之并. 注意到，定义中的瓷砖不一定连通. 对于连通的瓷砖，在 2 维中称为多联骨牌（Polyomino），在 3 维中称为多联立方体（Polycube）.

关于固定参数的平移密铺问题的不可判定性，一个显而易见的性质是：若问题对 $(n_0, k_0)$ 是不可判定的，则对任意 $n \geqslant n_0$ 和 $k \geqslant k_0$，问题对固定参数 $(n, k)$ 也是不可判定的.

文献中最早的关于固定瓷砖数目的平移密铺不可判定性结果由 Ollinger 在 2009 年对 2 维平面的情形给出. Ollinger 证明了 $(n, k) = (2, 11)$ 的平移密铺问题是不可判定的[9]. Ollinger 的证明方法是把王浩瓷砖的密铺问题归约为一组 11 个多联骨牌的密铺问题，其归约的核心框架和思想可用"总体有刚性、局部可选择"来概括. 由于平面中存在一组 8 个多联骨牌的非周期组[1]，Ollinger 在 2009 年提出 $(n, k) = (2, 8)$ 时的平移密铺问题也是不可判定的猜想[9]. 这一猜想被杨超和张诸俊的一系列工作证明了[14~16]. 杨超和张诸俊进一步证明 $(n, k) = (2, 7)$ 的不可判定性[20]. Kim 证明了 $(n, k) = (2, 5)$ 的不可判定性[7]. 本段前述的平面中的不可判定性结果都是在瓷砖为连通的多联骨牌的条件下得到的. 若允许瓷砖不连通，杨超和张诸俊证明 $(n, k) = (2, 4)$ 时平移密铺问题的不可判定性[21].

在维数 $n \geqslant 3$ 的情况下，杨超和张诸俊先后证明了 $(n, k) = (3, 6)$、$(n, k) = (4, 4)$ 和 $(n, k) = (4, 3)$ 的不可判定性[17~19]. 即目前关于固定维数的平移密铺问题的不可判定性的已知结果中，瓷砖最少数目为 $k = 3$. 本文的主要工作是把平移密铺问题的不可判定性的参数从 $(n, k) = (4, 3)$ 改进到 $(n, k) = (3, 3)$，即证明如下定理.

**定理 2** 平移密铺问题对参数 $(n, k) = (3, 3)$ 是不可判定的.

为了证明定理2，要把杨超和张诸俊在证明 $(n, k) = (2, 4)$ 一文[21]提出的非等距排列模拟王砖的方法推广应用到本文讨论的 3 维情形. 本文余下部分组织如下. 第2节介绍非等距排列需要用到的

关键组合引理. 第3节介绍构造多联立方体瓷砖所需的基本构件. 第4节在基本构件的基础上介绍 3 个多联立方体的具体构造方法. 第5节通过把任意一组王砖归约为一组 3 个多联立方体完成定理2的证明. 第6节是总结和展望.

## 2 组合引理

本文主要定理（定理2）的证明依赖于本节定义的不等距正整数集的存在性. 不等距集等价于 Golomb 尺（Colomb Ruler）和 Sidon 集，在组合和编码理论中被广泛地研究[3].

**定义 2** (不等距集) 设 $A \subset \mathbb{Z}^+$ 是一个有限集. 若对任意的 $i,j,i',j' \in A$, $i \neq j$, $i' \neq j'$, 且 $\{i,j\} \neq \{i',j'\}$, 都有 $|i-j| \neq |i'-j'|$, 则称 $A$ 为不等距集.

换言之，一个集合 $A$ 为不等距集，当且仅当 $A$ 中任意两对不同的整数对的距离不相等.

**性质 1** 不等距集存在.

**证明** 设 $n$ 为正整数，令 $B = \{1, 2^1, 2^2, \cdots, 2^{n-1}\}$，下证 $B$ 为不等距集. 任取 $0 \leqslant i, j, i', j' \leqslant n-1$, 且 $\{i,j\} \neq \{i',j'\}$. 不失一般性，设 $i < j$, $i' < j'$, $i' \leqslant i$. （反证）反设 $B$ 中两对正整数 $(2^i, 2^j)$ 与 $(2^{i'}, 2^{j'})$ 等距，即 $2^j - 2^i = 2^{j'} - 2^{i'}$, 则有

$$1 = 2^{j'-i'} - 2^{j-i'} + 2^{i-i'}.$$

若 $i = i'$, 上式化为 $1 = 2^{j'-i'} - 2^{j-i'} + 1$, 从而 $j = j'$, 矛盾. 若 $i \neq i'$, 上式左边为奇数，右边为偶数，亦矛盾. 证毕. □

上述证明中构造的不等距集的元素的最大距离 $\max(B) - \min(B) = 2^{n-1} - 1$ 相对于集合的大小 $n = |B|$ 呈指数关系. 定理3表明存在更小的不等距集.

**定理 3** (见 [11]) 存在 $n$ 元不等距集 $A = \{a_1, a_2, \cdots a_n\}$, 使得 $\max(A) - \min(A) = O(n^2)$.

定理2的证明需要用到比不等距集更强的模不等距集.

**定义 3** (模不等距集) 设 $A \subset \mathbb{Z}^+$ 是一个有限集，且 $m > \max(A)$. 若对任意的 $i,j,i',j' \in A$, $i \neq j$, $i' \neq j'$, 且 $\{i,j\} \neq \{i',j'\}$, 都有 $i - j \neq i' - j' \pmod{m}$, 则称 $A$ 为模 $m$ 不等距集.

模 $m$ 不等距集的条件可等价地表述为：对任意 $i < j$, $i' < j'$, 且 $\{i,j\} \neq \{i',j'\}$, 都有 $j - i \neq j' - i'$ 且 $j - i \neq m - (j' - i')$. 由此可见，模不等距集一定是不等距集，但反之不然.

**性质 2** 模不等距集存在.

**证明** 设 $n \geqslant 2$ 为正整数，令 $B = \{2^2, 2^3, \cdots, 2^n\}$, $m = 2^n + 2$, 下证 $B$ 为模 $m$ 不等距集. 任取 $2 \leqslant i, j, i', j' \leqslant n$, 且 $\{i,j\} \neq \{i',j'\}$. 不失一般性，设 $i < j$, $i' < j'$, $i' \leqslant i$. 由性质1的证明知：$2^j - 2^i \neq 2^{j'} - 2^{i'}$. 另外，亦有

$$2^j - 2^i \neq m - (2^{j'} - 2^{i'}).$$

这是因为上式左边能被 4 整除，右边为 $m - (2^{j'} - 2^{i'}) = (2^n + 2) - (2^{j'} - 2^{i'})$ 不能被 4 整除. □

**注释 1** 文 [21] 证明平移密铺问题在参数 $(n,k) = (2,4)$ 时的不可判定性，只需用到性质1不等距集的存在性. 本文证明定理2需要用到稍强的性质2模不等距集的存在性.

**注释 2** 本文在证明定理2时，没有特别优化所构造的 3 个多联立方体的大小. 若利用定理3及其模不等距集的版本，可以大幅度减小所需多联立方体的大小.

## 3 基本构件

若干个 $1\times 1\times 1$ 的单位正立方体面对面粘合起来形成的连通的立体图形称为多联立方体（Polycube）. 由 $n\times n\times n = n^3$ 个单位正立方体组成的边长为 $n$ 个单位的正立方体称为 $n$ 阶立方体. 本节将定义几类多联立方体作为基本构件. 这些基本构件将在下一节中被用来构造相对更加复杂的多联立方体.

多联立方体上若加上特殊形状的凹凸[1]，可确保某种刚性的密铺结构. 本文主要采用的凹凸形状为图2所示的"十"字形多联立方体. 此"十"字形多联立方体也是定理2的证明中所需的三个多联立方体之一，称之为填充块. 填充块由 9 个单位正立方体在同一水平高度中组合成"十"字形，即图2为俯视图. 注意到，基本凹凸形状有无限多种不同的选择，此处我们选取与 Kim 文 [7] 一致的相对比较简单的"十"字形.

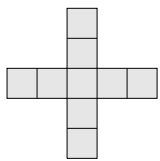

**图 2 填充块**
**Figure 2** A Filler

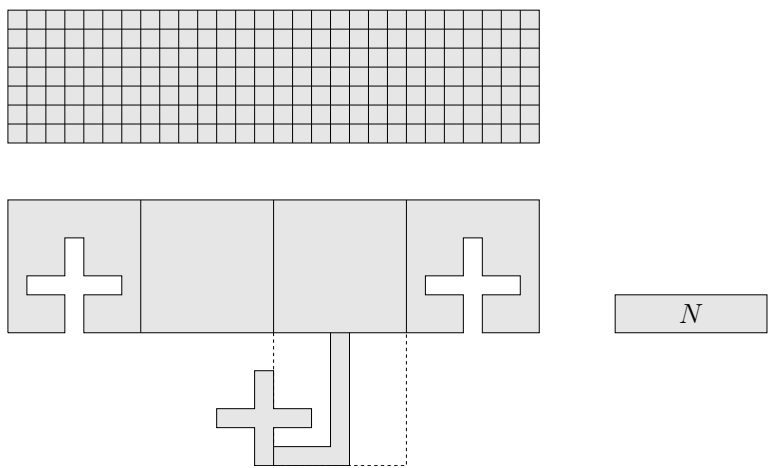

**图 3 构件 $N$ 的一级分层图（左）及其示意图（右）**
**Figure 3** The level-1 layer diagram of the building block $N$ (left), and its symbolic representation (right)

第一类基本构件是在 $28\times 7\times 7$（即长、宽、高分别是 28、7、7，亦可视为东西方向排成一行的 4 个 7 阶立方体）的多联立方体的基础上添加若干凹凸得到. 构件 $N$ 由图3的分层图定义. 构

---
1) 类似于中国传统木建筑中的榫卯结构.

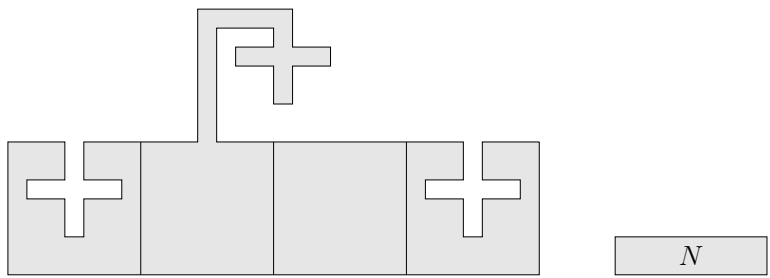

**图 4  构件 $N$ 的另一朝向（左）及其示意图（右）**
**Figure 4** Another orientation of building block $N$ (left), and its symbolic representation (right)

件 $N$ 的第 1、2、3、5、6、7 层均为 $28 \times 7 \times 1$ 的多联立方体，如图3（左上）所示. 凹凸结构只加在构件 $N$ 的最中间的一层，即第 4 层，如图3（左下）所示：左右两侧各加了一个"十"字形的凹进形状，中间加了一个 $L$ 字形和一个"十"字形结合的凸出形状. $L$ 字形的两边长度分别为 7 和 4. 构件 $N$ 整体地用图3（右）的示意图表示，将在下一节中用于构造更大的多联立方体.

构件 $N$ 的凹凸结构除了可以加在南向（如图3所示），也可以加在北向（如图4）. 两种朝向的构件 $N$ 的形状是完全一样的，一种朝向旋转 180° 后和另一种朝向重合. 本节后面定义的基本构件 $F$, $E$, $M$ 和 $M^+$ 同样也有两种朝向. 为了节省篇幅，只图示其中一种朝向，旋转 180° 便得到另外一种朝向. 此外，本节后面定义的基本构件 $F$, $E$, $M$ 和 $M^+$ 的第 1、2、3、5、6、7 层和构件 $N$ 完全一样，都是 $28 \times 7 \times 1$ 的多联立方体. 不同之处只在于中间的第 4 层，故对这几个构件的图示只画出第 4 层.

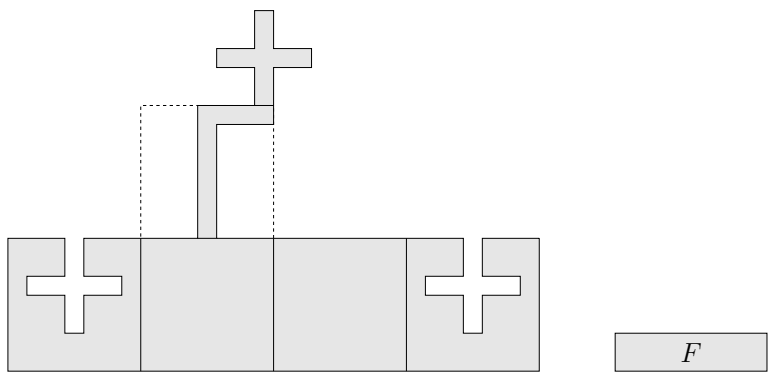

**图 5  构件 $F$（左）及其示意图（右）**
**Figure 5** The building block $F$ (left) and its symbolic representation (right)

构件 $F$ 的第 4 层如图5（左）所示. 构件 $F$ 两侧的"十"字形凹进和构件 $N$ 相同. 构件 $F$ 与构件 $N$ 唯一不同之处是中间的"十"字形的凸起附加在 $L$ 字形凸起上的位置不一样. 在构件 $N$ 中，"十"字形凸起的附加位置距离构件的主体部分更近（Near），而构件 $F$ 中，附加的位置更远（Far）.

构件 $E$ 的第 4 层如图6（左）所示，它可视为构件 $N$ 或构件 $F$ 去掉"十"字形的凹进或凸起，只保留 $L$ 字形的凸起.

图7是两个背靠背的构件 $M$ 的中间第 4 层. 构件 $M$ 在中间层的南向或北向添加一个 $L$ 字形的

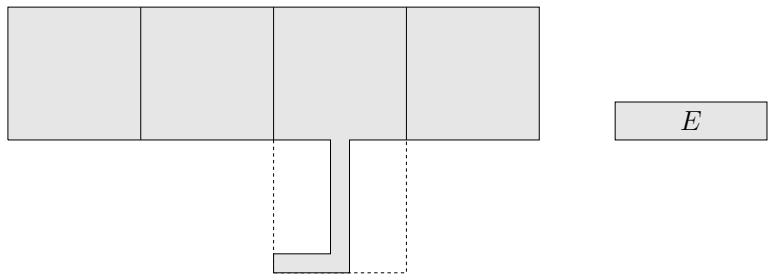

**图 6** 构件 $E$（左）及其示意图（右）
**Figure 6** The building block $E$ (left) and its symbolic representation (right)

凹进和一个"十"字形的凹进的结合体. 两个构件 $M$ 背靠背排列之后, 凹进的部分连成一个整体, 使得凹进部分要么可以刚好被两个构件 $N$ 的凸起部分填满, 要么被两个构件 $F$ 的凸起部分填满.

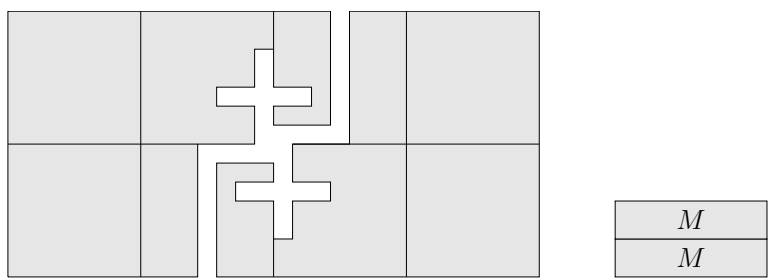

**图 7** 两个构件 $M$（左）及它们的示意图（右）
**Figure 7** Two building blocks $M$ (left) and their symbolic representations (right)

在构件 $M$ 的基础上, 若再在两侧加上两个"十"字形凸起, 就得到了构件 $M^+$. 两个背靠背的构件 $M^+$ 的中间第 4 层如图8（左）所示.

容易验证, 以上定义的基于 4 个 7 阶立方体基本构件都是连通的. 对于构件 $M$ 和 $M^+$, 虽然它们的第 4 层单独看不连通, 但是通过第 1、2、3、5、6、7 层粘合起来后整体是连通的.

此外, 前面定义的基本构件 $N$, $F$, $E$, $M$ 和 $M^+$ 都有向北或向南两种朝向, 然而在简化的示意图中省略了朝向的标记. 在下一节用基本构件来构造更大的多联立方体时, 每一个基本构件的南侧和北侧只有一侧作为它所组成的更大的多联立方体的外侧. 因此, 默认所有的基本构件都朝向向外的一侧 (即凹凸结构是添加到向外的一侧), 即使省略朝向的标记也不会引起歧义. 如图8（右）两个基本构件 $M^+$ 组成的多联立方体中, 北边的构件朝北, 南边的构件朝南.

除了第一类基于 4 个 7 阶正立方体的基本构件之外, 我们还定义第二类共 3 对基于 1 个 7 阶正立方体的基本构件. 第一对 $X^+$ 和 $X^-$ 如图9所示, 它们分别是在 7 阶正立方体的基础上添加一个凸起和一个凹进. 构件 $X^+$ 和 $X^-$ 的 1、2、3、5、6、7 层都是 $7 \times 7$ 的 49 个单位正立方体（如图9左）. 构件 $X^+$ 的第 4 层向东方向 (即 $X$ 轴的正方向) 有一个凸起 (如图9中上). 而 $X^-$ 的第 4 层的西侧有一个凹进 (如图9中下), 其形状和 $X^+$ 第 4 层东侧的凸起完全一样. 和前面第一类构件的"十"字形凹凸形状类似, 此处 $X^+$ 和 $X^-$ 的凹凸形状的选取也有很大的自由度, 只要显著地不同于第一类构件中的"十"字形和 L 字形即可. 不同的凹凸形状实现不同的功能模块. 若把 $X^+$ 和 $X^-$ 绕某垂直旋转轴旋转 90°, 使得凸起在北侧 (即 $Y$ 轴的正方向), 凹进在南侧, 则把得到的

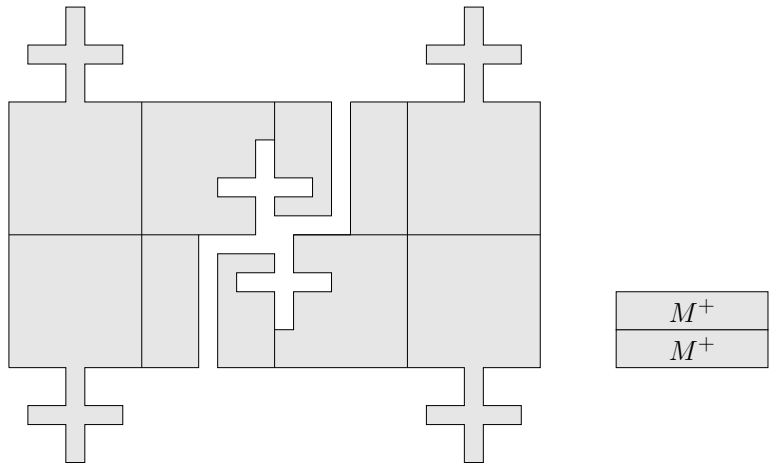

**图 8** 两个构件 $M^+$（左）及它们的示意图（右）
**Figure 8** Two building blocks $M^+$ (left) and their symbolic representations (right)

基本构件分别称为 $Y^+$ 和 $Y^-$. 若把 $X^+$ 和 $X^-$ 绕某平行于 $Y$ 轴的直线旋转 90°，使得凸起在上侧（即 $Z$ 轴的正方向），凹进在下侧，则把得到的基本构件分别称为 $Z^+$ 和 $Z^-$.

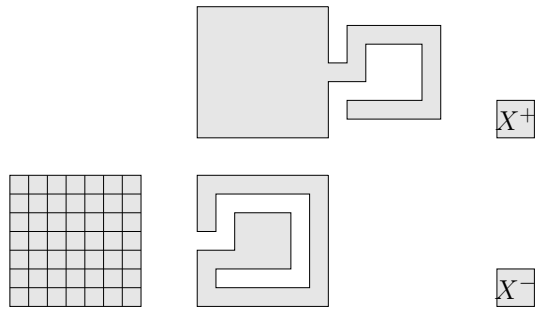

**图 9** 基本构件构件 $X^+$ 和 $X^-$（左、中）及它们的示意图（右）
**Figure 9** Building blocks $X^+$, $X^-$ (left and middle) and their symbolic representations (right)

## 4 一组 3 个多联立方体

本节以图1的一组 3 个王砖为例，用第3节所定义的基本构件来构造与该王砖组对应的一组 3 个多联立方体瓷砖：填充块、编码块、连接块. 构造的方法可自然地应用到一般的一组 $n$ 个王砖的情形.

第一个多联立方体填充块（见图2）已经在第3节介绍了. 填充块的大小是固定的，不随王砖组的瓷砖数目或颜色数目的变化而变化.

第二个多联立方体编码块通过某种有规律的方式，把一组 $n$ 个王砖都编码到一个多联立方体中. 因此，编码块的大小随着王砖组的瓷砖数目和颜色数目的增大而增大. 编码块具有多层结构，其层数较大，由有限的几类分层结构通过特定方式在垂直方向排列组合而成. 对应图1所示王砖组的编码块的分层结构如图10所示. 图10展示了编码块的其中四层，按图10中从上到下的顺序分别为：编码层、结构层、非编码层、结构层. 其中编码层和非编码层统称功能层. 需要指出的是，本节的分层示

意图（图10编码块和图11连接块）的每一层的高度为 7 个单位, 即基本构件 $N$、$F$、$E$ 等的高度. 而第3节中基本构件 $N$、$F$ 和 $E$ 等的分层图中, 每一层的高度为 1, 即一个 $1 \times 1 \times 1$ 单位正立方体的高度. 为了区别两种不同层次的分层图, 第3节的分层图称为一级分层图（如图3）, 本节的分层图称为二级分层图.

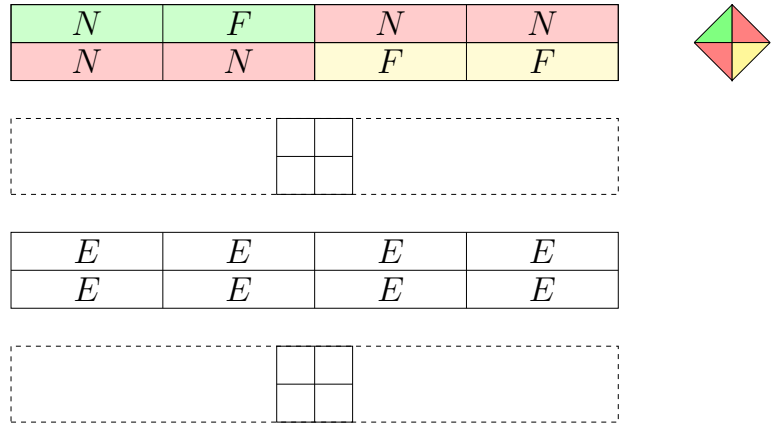

**图 10 编码块的二级分层示意图**
**Figure 10** The level-2 layer diagram of encoder

对应图1所示王砖组的编码块一共有 $1026 = 2^{10} + 2$ 层, 按空间位置从下到上依次编号为 1 至 1026. 其中奇数层为结构层, 偶数层为功能层, 即结构层和功能层在垂直方向交错排列. 在功能层中, 只有编号为 $2^k (2 \leqslant k \leqslant 10)$ 的为编码层, 其余均为非编码层. 每个编码层对一个王砖进行编码, 其中第 4、8、16 层编码图1的第一个王砖, 第 32、64、128 层编码图1的第二个王砖, 第 256、512、1024 层编码图1的第三个王砖. 换言之, 每个王砖在编码块中被重复编码三次. 注意到, 编码层所在层数之集 $\{2^2, 2^3, \cdots, 2^{10}\}$ 相对于总层数 $2^{10} + 2$ 是一个模不等距集. 结构层均为 $14 \times 14 \times 7$ 的长方体, 即由 4 个 7 阶正方体组成, 居中与其余各层对齐. 非编码层由 2 行 4 列共 8 个基本构件 $E$ 组成. 编码层则由 2 行 4 列的共 8 个基本构件 $N$ 或 $F$ 组成. 编码层的 8 基本构件分层 4 组, 每组为东西方向相邻的 2 个基本构件, 分别编码了一个王砖的一条边. 图10中的编码层编码了图1的第一个王砖. 因为图1所示王砖组一共有 4 种不同的颜色, 因此两个东西方向排列的 2 基本构件 $N$ 或 $F$ 刚好能以二进制的方式对 4 种不同的颜色编码. 连续的两个基本构件 $NN$、$NF$、$FN$ 和 $FF$ 分别编码颜色红、绿、蓝和黄.

下面描述一般情况下, 一组王砖所对应的编码块的结构. 设有一组 $n$ 个王砖, 其边共有 $m$ 种不同的颜色, 并令 $t = \lceil \log_2 m \rceil$（称此王砖组以 $n, m, t$ 为参数）. 则该王砖组对应的编码块有 $2^{3n+1} + 2$ 层, 同样地奇数层为结构层, 偶数层为功能层. 功能层由 2 行 $2t$ 列共 $4t$ 个基本构件组成, 其中非编码层的 $4t$ 个基本构件均为 $E$, 而编码层的 $4t$ 个基本构件为 $N$ 或 $F$. 类似于图10, 每个编码层的 $4t$ 个基本构件分为四组, 每组为东西方向连续的 $t$ 个基本构件 $N$ 或 $F$, 且每组编码一个王砖的一条边（即 $t$ 位的 $N, F$ 序列足以编码 $2^t \geqslant m$ 种不同颜色）, 四组合起来编码一个王砖. 对 $1 \leqslant i \leqslant n$, 第 $2^{3i-1}$, $2^{3i}$ 和 $2^{3i+1}$ 层编码第 $i$ 个王砖. 即编码层数之集构成一个模 $2^{3n+1} + 2$ 不等距集. 其余编号非 $2^k (2 \leqslant k \leqslant 3n+1)$ 的功能层（即偶数层）为非编码层. 一般情况的结构层（即奇数层）和图10中的结构层完全一样, 由 4 个 7 阶正方体组成, 位于该层的正中.

第三个多联立方体称为连接块,主要起到固定密铺的整体框架,并且确保三个多联立方体的密铺模拟了所对应的王砖组的平面密铺的作用. 和编码块一样,连接块也是多层结构,由较少的几类分层结构按一定方式垂直排列而成. 对应图1的王砖组的连接块的其中若干分层结构如图11所示. 图11中从上到下所展示的四层分别为匹配层、结构层、非匹配层、结构层. 其中匹配层和非匹配层统称功能层. 匹配层由 2 行 4 列共 8 个基本构件 $M^+$ 组成,非匹配层由 2 行 4 列共 8 个基本构件 $M$ 组成. 结构层在 2 行 16 列共 32 个 7 阶立方体的基础上,北侧居中位置再排列 2 行 14 列共 28 个 7 阶立方体. 整个连接块一共有 1026 层,和相应的编码块具有相同的层数. 和编码块一样,把连接块按照空间位置从下到上对各层依次编号为 1 至 1026,且奇数层为结构层,偶数层为功能层. 在功能层中,第 1026 层为匹配层,其余的功能层均为非匹配层. 此外,连接块还有 6 个 7 阶立方体被特殊的构件替代(见图11):第 1 层的第 2 行的第 1 列和第 16 列分别为构件 $X^-$ 和构件 $X^+$,第 1 层的第 1 行的第 1 列为构件 $Y^-$,第 1 层的第 4 行的第 9 列为构件 $Y^+$,第 1 层的第 1 行的第 16 列为构件 $Z^-$,第 1026 层的第 1 行的第 16 列为构件 $Z^+$. 其中构件 $Z^+$ 是叠加在原来的构件 $M^+$ 上,原构件 $M^+$ 在此位置有南侧的一个"十"字形凸起,构件 $Z^+$ 则在此位置添加了上侧的一个凸起,两个凸起并行不悖.

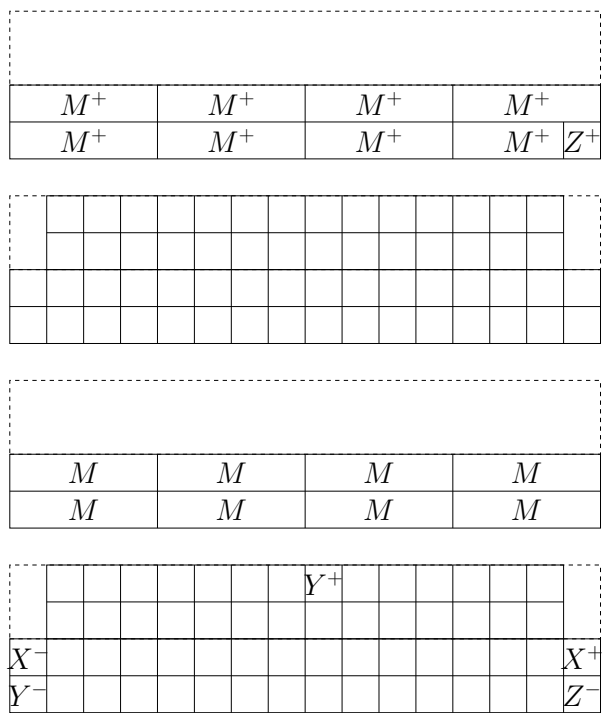

**图 11  连接块的二级分层示意图**
**Figure 11**  The level-2 layer diagram of linker

一般情况下,对一组参数为 $m$, $n$, $t$ 的王砖,其对应的连接块层数为 $2^{3n+1}+2$. 其中奇数层为结构层,偶数层为功能层. 第 $2^{3n+1}+2$ 层为匹配层,其余的功能层均为非匹配层. 匹配层由 2 行 $2t$ 列共 $4t$ 个基本构件 $M^+$ 组成,非匹配层由 2 行 $2t$ 列共 $4t$ 个基本构件 $M$ 组成. 结构层由主体部分和延申部分两部分组成. 主体部分是 2 行 $8t$ 共 $16t$ 个 7 阶正立方体构成的 $56t \times 14 \times 7$ 的长方体.

主体部分在空间垂直方向与功能层对齐. 延申部分是在主体部分北侧的 2 行 $8t-2$ 列共 $16t-4$ 个 7 阶正立方体. 延申部分相对主体部分居中对齐, 即东西两侧各缩进 1 个 7 阶正立方体. 连接块的 6 基于 7 阶正立方体的特殊的构件的位置为: 第 1 层的第 2 行的第 1 列和第 $8t$ 列分别为构件 $X^-$ 和构件 $X^+$, 第 1 层的第 1 行的第 1 列为构件 $Y^-$, 第 1 层的第 4 行的第 $4t+1$ 列为构件 $Y^+$, 第 1 层的第 1 行的第 $8t$ 列为构件 $Z^-$, 第 $2^{3n+1}+2$ 层的第 1 行的第 $8t$ 列为构件 $Z^+$.

综上所述, 对任意的一组王砖, 我们给出一个系统的方法构造一组与之相对应的 3 个多联立方体瓷砖: 填充块、编码块和连接块. 其中填充块大小固定不变, 编码块和连接块的大小随着王砖组的参数 $n$ 呈指数增长.

## 5 定理2的证明

**证明** 在上一节中, 我们已经对任意一组王砖 $W$, 相应地构造了一组 3 个多联立方体 $P$: 填充块、编码块和连接块. 以下只需证明, 王砖组 $W$ 能平移密铺平面当且仅当相应的多联立方体组 $P$ 能平移密铺 3 维空间.

- 用多联立方体组 $P$ 平移密铺 3 空间时, 连接块必须使用. 在平面中, 可平移密铺整个平面的单一瓷砖有充分必要的刻画条件[4]. 特别地, 对于单个多联骨牌能否平移密铺平面, 有快速算法[13]. 由平面单个瓷砖的平移密铺条件, 易知图2 "十" 字形填充块的水平面投影无法平移密铺平面. 进而, 只用填充块无法密铺整个空间 (否则, 可用一水平截面得到 "十" 字形的平面密铺, 矛盾). 因此编码块或连接块必须使用. 若使用了编码块, 则编码块上的基本构件 $N, F$ 上的 $L$ 字形和 "十" 字形的组合体凸出部分必须要用连接块的相同形状的凹进匹配, 因此连接块也必须使用.

- 连接块在密铺时形成具有格 (Lattice) 结构的刚性框架. 由于一个连接块第 1 层中的构件 $X^+$ 只能和另一个连接第 1 层中的构件 $X^-$ 对接, 才能在构件 $X^+$ 局部不留下空隙, 所以这两个连接块在水平向是完全对齐的, 即一个连接块的第 $i$ 层对齐另一个连接块的第 $i$ 层, 对每个 $1 \leqslant i \leqslant 2^{3n+1}+2$. 又因为一个连接块第 1 层中的构件 $Y^+$ 只能和另一个连接块第 1 层中的构件 $Y^-$ 对接, 故南北方向相邻的两个连接块有东西方向的相对位移. 由 $Y^+$ 和 $Y^-$ 构件的位置, 东西方向的相对位移刚好是连接块的东西方向长度的一半. 因此在密铺的水平截面中, 连接层的排列如图12 (结构层) 或图13 (功能层) 中灰色块所示. 再由一个连接层第 $2^{3n+1}+2$ 层的构件 $Z^+$ 只能和另一个连接第 1 层的构件 $Z^-$ 对接, 故连接层的排列模式向下或向上无限延申都是完全相同的. 若以一个 7 阶的正立方体为一个基本单元, 连接块的如上所述的排列方式构成一个以 $(16t, 0, 0)$, $(2, 8t, 0)$ 和 $(0, 0, 2^{3n+1}+2)$ 为一组基的格. 在此格结构中, 连接块之间的相对位置固定无法移动, 这就是第1节引言中提到的所谓的 "总体有刚性".

- 编码块放置在连接块形成的格结构框架下的空间中, 形成垂直方向双向无限的柱状序列. 由于连接块中的 $L$ 字形凹进只能用编码块中的 $L$ 字形凸起才能填满, 所以在连接块形成的格结构框架中, 剩下的空间必须用编码块填充, 如图12和图13中的紫色部分 (含绿色部分). 并且, 在垂直方向上的, 一连串编码块完全对齐, 首尾相连, 中间不留空层, 形成上下双向无限的一个序列. 这是因为若两个垂直对齐的编码块中间有空层, 当空层的层数大于或等于一个编码块的层数时, 由 $L$ 字形凹凸的相互耦合, 空层只能再填放另一个编码块; 当空层的层数小于一个编码块的层数时, 空层则无法填满了. 如图13所示, 每个紫色的矩形区域对应一个垂直方向双向无限的编码块序列.

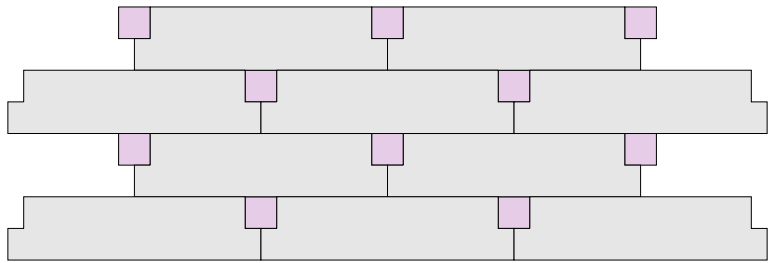

**图 12** 密铺的结构层的水平截面

**Figure 12** The horizontal sectional plane of structural layers in a tiling

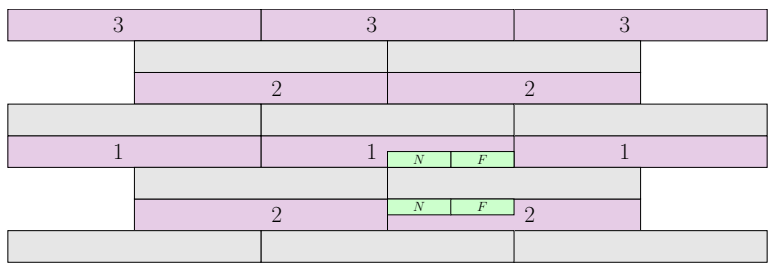

**图 13** 密铺的功能层的水平截面

**Figure 13** The horizontal sectional plane of functional layers in a tiling

• 每串由编码块形成的柱状序列相对连接块的格框架，在垂直方向可独立地上下调整位置，选择某一个编码层和连接块唯一的匹配层对齐. 当我们只考虑其中一连串垂直方向双向无限的编码块序列时，由于连接块匹配层存在"十"字形凸起，所以必须和编码块中的编码层对齐. 而非匹配层既可以和编码层对齐，也可以和非编码层对齐，无论哪一种情况，都不会产生重叠（但可能会留下小空隙）. 而编码块中有多个编码层，故单独地看一连串编码块序列时，可以对整个序列上下调整位置，选择哪一个编码层和连接块中唯一的匹配层对齐. 这就是第1节引言中提到的所谓的"局部可选择".

• 所有的匹配层的密铺模式都是完全相同的，模拟了同一种王砖的密铺方案. 由于编码块和连接块的高度（或层数）是相等的，所以当一个编码块的某一个编码层与匹配层对齐时，此编码块所属的垂直方向一连串编码块序列都是同一个编码层与匹配层对齐. 因此在整个空间密铺中，每一个匹配层的水平截面都是完全相同的. 在匹配层的水平截面中，每一个编码块的编码层与四个其它编码块的编码层相邻，分别位于东北、东南、西南、西北方向. 这模拟了一个王砖和四个王砖相邻. 同时，相邻的编码层不会产生重叠，当且仅当相邻的基本构件是 $N$ 与 $N$ 对应，且 $F$ 与 $F$ 对应. 这模拟了相邻的王砖的相邻边的颜色要一致. 如图13中绿色部分就示意了编码层中两个相邻的模拟王砖的相邻边的颜色相同均为绿色.

• 非匹配层的密铺模式可实现不重叠. 在匹配层可以实现密铺无重叠的条件下，本段论证每个非匹配层也可以实现无重叠. 由第4节编码块的构造可知，编码层所在层数的全体，即 $\{2^2,\cdots,2^{3n+1}\}$，构成一个模 $2^{3n+1}+2$ 不等距集. 因此，相邻的两个编码块，要么编码层完全对齐（即整个编码块对齐，第 1 层对齐第 1 层，第 2 层对齐第 2 层，等等），要么仅有一个编码层对齐（即与连接块的匹配层对齐的编码块）. 只需证明可以避免相邻两个编码块完全对齐的情况出现即可. 按以下顺序确定每一个编码块相对于连接块的位置（即相对高度位置）. 先确定东西方向的一行编码块的位置. 因

为东西方向位于同一行的编码块（从模拟王砖密铺的角度看）是不相邻的，位置可以任意确定. 然后向北和向南逐一确定其它每个编码块序列的相对高度位置. 图13中标示的 1, 2, 3 即表明了编码块序列的位置确定的顺序：从中间向南北延伸. 在这一过程中，每个待确定位置的编码块与两个已确定位置的编码块相邻，即只有 2 种情况会和这两个相邻的编码块完全对齐. 又注意到，每个王砖在编码块中被编码了 3 次，即待确定位置的编码块有 3 种不同的高度位置选择，使得与匹配层对齐的编码层都是编码了同一个王砖. 由抽屉原理，可在不改变（待确定位置的编码块序列）与匹配层对齐的编码层所模拟的王砖的条件下，避免整个编码块完全对齐. 若单独地从某一个非匹配层的水平截面看，就如一个个孤立的编码层被包围在非编码层的汪洋大海中，不存在两个编码层相邻的情况. 从而，在非编码层也不会出现一般编码层基本构件 $N$ 和另一编码层的基本构件 $F$ 相邻的情况，所以不会产生重叠.

- 剩下的小空隙都可以用"十"字形填充块填满. 经过以上的步骤之后，剩下只有"十"字形的空隙，都可以用填充块填满.

由以上论证可知，这一组 3 个多联立方体能密铺整个空间当且仅当匹配层无重叠. 而匹配层模拟了王砖的平面密铺，从而又等价于这组 3 个多联立方体所模拟的王砖组能密铺平面. 由王砖密铺问题的不可判定性（定理1）可得，一组 3 个多联立方体平移密铺整个 3 维空间是不可判定的. □

## 6 总结与展望

本文证明了一组 3 个多联立方体平移密铺 3 维空间的不可判定性，这是目前已知关于固定参数 $(n, k)$ 的平移密铺问题的不可判定性中 $k$ 最小的结果. 这一结果的证明关键是利用了模不等距集的存在性，在编码块上不等距地排列编码层. 注意到文 [19] 证明平移密铺问题在参数 $(n, k) = (4, 3)$ 时的不可判定性，其三个块对应地也可以称为编码块、连接块和填充块. 文 [19] 是利用了第 4 维空间的错位来实现只需一个连接块. 而本文创新地采用模不等距排列，在更低的 3 维空间就就实现了只需一个连接块. 模不等距排列的方法也许还有潜力，若能更有效地与其它方法结合，有可能得到固定维数空间 $k = 2$ 个瓷砖的平移密铺问题的不可判定性.

# On the Undecidability of Tiling the 3-dimensional Space with a Set of 3 Polycubes

Chao YANG[1*] & Zhujun ZHANG[2]

1. *School of Mathematics and Statistics, Guangdong University of Foreign Studies, Guangzhou* 510006, *China*;
2. *Big Data Center of Fengxian District, Shanghai* 201499, *China*
* Corresponding author. E-mail: yangchao@gdufs.edu.cn

**Abstract**  Translational tiling problems are among the most fundamental and representative undecidable problems in all fields of mathematics. Greenfeld and Tao obtained two remarkable results on the undecidability of translational tiling in recent years. One is the existence of an aperiodic monotile in a space of sufficiently large dimension. The other is the undecidability of translational tiling of periodic subsets of space with a single tile, provided that the dimension of the space is part of the input. These two results support the following conjecture: there is a fixed dimension $n$ such that translational tiling with a single tile is undecidable. One strategy towards solving this conjecture is to prove the undecidability of translational tiling of a fixed dimension space with a set of $k$ tiles, for a positive integer $k$ as small as possible. In this paper, it is shown that translational tiling the 3-dimensional space with a set of 3 polycubes is undecidable.

**Keywords**  undecidability, translational tiling, polycube, Wang tile, Golomb ruler

**MSC(2020)**  52C22, 68Q17